\documentclass[4paper,12pt]{article}
\usepackage[left=2cm,right=2cm,top=2.5cm,bottom=2.5cm]{geometry}
\usepackage{graphicx, float} 
\usepackage{subfiles} 
\usepackage{amsmath, esint, amsthm, physics, amsfonts, amssymb} 

\usepackage[colorlinks=true, allcolors=blue]{hyperref} 

\usepackage[numbers,sort&compress]{natbib}
\usepackage{url}

\usepackage{indentfirst}
\usepackage{sectsty}
\usepackage{stmaryrd}
\usepackage{subfig}
\usepackage{lscape}
\usepackage{rotating}
\usepackage[stable]{footmisc} 
\usepackage{titling}
\usepackage{lipsum}
\usepackage{ragged2e}
\usepackage{enumitem}
\usepackage{authblk} 
\usepackage{orcidlink}
\usepackage{multicol, fancyhdr, xcolor} 
\usepackage{tikz}

\theoremstyle{definition}

\theoremstyle{remark}

\theoremstyle{proposition}
\newtheorem{proposition}{Proposition}[section]

\DeclareMathAlphabet\mathbfcal{OMS}{cmsy}{b}{n}

\title{Carleman Linearization of Partial Differential Equations}
\author{Tamás Á. Vaszary\footnote{Email: \href{mailto:Tamas.Vaszary@maths.ox.ac.uk}{Tamas.Vaszary@maths.ox.ac.uk}}  \orcidlink{0000-0001-7811-222X}}
 \affil{Mathematical Institute, University of Oxford, 
United Kingdom}
\date{November 2024}

\begin{document}

\maketitle

\begin{abstract}
    Carleman linearization is a technique that embeds systems of ordinary differential equations with polynomial nonlinearities into infinite dimensional linear systems in a procedural way. In this paper we generalize the method for systems of partial differential equations with quadratic nonlinearities, while maintaining the original structure of Carleman linearization. Furthermore, we apply our approach to Burger's equation and to the Vlasov equation as examples.
\end{abstract}

\section{Introduction}
Finding approximate global linear representations of nonlinear dynamical systems has become increasingly important for computational purposes. This is especially relevant for the emerging sector of quantum computation, which has been proved to have a performance advantage over classical methods for numerous linear problems \cite{hhl,an2023theoryquantumdifferentialequation,coppersmith2002approximatefouriertransformuseful}.

For systems of ordinary differential equations (ODEs), Koopman operator theory provides a general framework for embedding the dynamics into a possibly infinite dimensional linear system \cite{koopman0,koopman1, koopman2, koopman3}. That system evolves a vector of so-called measurement functions of the original variables. The drawback of this procedure is that finding suitable measurement functions and the corresponding linear operator cannot be carried out for arbitrary systems and it is not procedurally applicable.

A special case of Koopman theory is Carleman linearization, which restricts the measurement functions to be monomials of the original dynamical variables that arise from Kronecker powers of the state vector (i.e. Kronecker products of the state vector with itself) \cite{carleman0,carleman1,forets2017explicit,Arash}. This way the linear operator can be expressed analytically and the system becomes procedurally implementable. As a result, it has already been used in real world engineering problems \cite{carleman3} and quantum algorithms \cite{liu,liu2,surana2023efficient,krovi}.

Generalizing this framework for systems of partial differential equations (PDEs) would carry key importance, due to the resulting linear PDE system potentially becoming solvable on analog quantum simulators. This would eliminate coordinate discretization error completely and would potentially offer a computational advantage as well in terms of complexity.

One piece of previous work already identified the key idea required for this generalization, namely constructing continuous Kronecker powers of the state vector, using different copies of the independent variables in each term \cite[Sections 3.5 \& 4.2]{Kowalski}. However, they do not provide a concrete, procedural way to implement this framework for general systems. The only other attempt at Carleman linearization of PDE systems was carried out without using the continuous analog of the Kronecker product mentioned above \cite{carleman_weird}. As such, this approach also lacks practical applicability.

In this paper we present a framework for the Carleman linearization of PDE systems, which has a similar structure to that of regular Carleman linearization and importantly, it also maintains the property of being procedurally implementable on a computer. We restrict ourselves to quadratic nonlinearities in this work, but it is known that higher dimensional polynomial nonlinearities can be transformed into quadratic ones within the regular Carleman framework \cite[Section 3.3]{forets2017explicit}. This is also expected to be true for the PDE case as discussed later.

For the sake of completeness, we review Carleman linearization in the next section, following the approach and notation of \cite{forets2017explicit} where applicable. Then, in Section \ref{sec:Carleman_for_pdes} we present our main results. Next, in Section \ref{sec:toy_problems} we apply the presented framework to toy problems as a proof of concept, namely to Burger's equation and the Vlasov equation. Finally, in Section \ref{sec:disc_conlc} we discuss the prospects and conclusion of this work.

\section{Review of Carleman linearization of ODEs}\label{sec:review}
\subsection{Formulation}\label{sec:formulation}
General ODE systems with quadratic nonlinearities can be written in the form
\begin{align}
\frac{{\rm d}\mathbf{u}}{{\rm d}t}=F_0+F_1\mathbf{u}+F_2\mathbf{u}^{\otimes 2}, \quad\quad \mathbf{u}(0)=\mathbf{u}_0,\label{eq:dudt_1}
\end{align}
where $\mathbf{u}=[u_1,\dots ,u_n]^\mathsf{T}\in\mathbb{R}^n$ is the state vector with elements $u_k(t)$ to be evolved over time $t\in[0,T]$, and where $\mathbf{u}^{\otimes2}:= \mathbf{u}\otimes\mathbf{u}=[u_1^2,u_1u_2,\dots,u_1u_n,u_2u_1,\dots,u_n u_{n-1},u_n^2]^\mathsf{T}\in \mathbb{R}^{n^2}$ contains the quadratic nonlinearities. The latter is the second Kronecker power of $\textbf{u}$. Furthermore $F_j \in\mathbb{R}^{n\times n^j}$ with $j=0,\:1,\:2$ are matrices that are assumed to be time independent.

To work within the Carleman linearization framework, we define the auxiliary variables $\mathbf{y}_i:=\mathbf{u}^{\otimes i}\in \mathbb{R}^{n^i}$ to be the $i$-th Kronecker power of $\mathbf{u}$. With this definition, $\mathbf{y}_i$ contains all possible monomials of the components of $\mathbf{u}$, such that the sum of the powers in the monomial is $i$. We now describe the Carleman linearization of ODEs of the form of \eqref{eq:dudt_1}. Note that we do not explicitly write out the time dependence of state variables from now for readability.

\begin{proposition}
The vectors $\mathbf{y}_i$ are evolved by the tridiagonal linear system\footnote{Note that for level $i=1$, the $j=0$ contribution on the RHS of \eqref{eq:dxidt_1} is inhomogeneous and is simply $F_0$. This means that with a slight abuse of notation $\mathbf{y}_0=\mathbf{u}^{\otimes 0}=1$ is understood, which is not a dynamical variable.}
\begin{align}
\frac{{\rm d}\mathbf{y}_i}{{\rm d}t}=\sum_{j=0}^2 A_{i+j-1}^i \mathbf{y}_{i+j-1},\label{eq:dxidt_1}
\end{align}  
where the composite matrices $A_{i+j-1}^i \in\mathbb{R}^{n^i\times n^{i+j-1}}$ are given by
\begin{align}
A_{i+j-1}^i:=\sum_{\nu=1}^i \overbrace{\mathbb{I}_{n} \otimes \cdots \otimes \underbrace{F_j}_{\nu\text{-th position }} \otimes \cdots \otimes \mathbb{I}_{n }}^{i \text { factors }},
\end{align}
with $\mathbb{I}_{n }$ denoting an $n\times n$ identity matrix.

\begin{proof}
To realize the form of the matrices above, we start by differentiating $\mathbf{y}_i$ and then use the Leibniz rule, obtaining
\begin{align}
\begin{split}
\frac{{\rm d}\mathbf{y}_i}{{\rm d}t}&=\frac{{\rm d}\mathbf{u}}{{\rm d}t}\otimes\mathbf{u}\otimes \dots \otimes \mathbf{u}+\dots +\mathbf{u}\otimes\dots\otimes \mathbf{u}\otimes \frac{{\rm d}\mathbf{u}}{{\rm d}t}\\
&= \sum_{\nu=1}^i \left(\mathbf{u}\otimes\dots\otimes\sum_{j=0}^{2}F_j \mathbf{u}^{\otimes j}\otimes\dots\otimes \mathbf{u}\right).
\end{split}
\end{align}
We can now factor out the operators acting on $\mathbf{y}_{i-1},\mathbf{y}_i$ and $\mathbf{y}_{i+1}$ as
\begin{align}
\begin{split}
\frac{{\rm d}\mathbf{y}_i}{{\rm d}t}&= \sum_{j=0}^2\sum_{\nu=1}^i
\left(\mathbf{u}\otimes\dots\otimes F_j\mathbf{u}^{\otimes j}\otimes\dots \otimes\mathbf{u}\right)\\
&=\sum_{j=0}^2\sum_{\nu=1}^i\left(\mathbb{I}_{n }\otimes \dots \otimes F_j\otimes\dots\otimes \mathbb{I}_{n }\right)\left(
\mathbf{u}\otimes\dots\otimes \mathbf{u}^{\otimes j}\otimes\dots\otimes \mathbf{u}
\right)\\
&=\sum_{j=0}^2 A_{i+j-1}^i \mathbf{y}_{i+j-1},
\end{split}
\end{align}
hence confirming \eqref{eq:dxidt_1}.
\end{proof}
\end{proposition}

\subsection{Truncation}\label{sec:truncation_ODE}
For the linearized system in \eqref{eq:dxidt_1} to be actually solvable on a computer, the system must be truncated at level $N$. Hence we define the Carleman linearized and truncated state vector $\mathbf{z}:=[\mathbf{z}_1^\mathsf{T},\dots , \mathbf{z}_N^\mathsf{T}]^\mathsf{T}\in \mathbb{R}^{\Delta}$ that has block elements $\mathbf{z}_m:=\mathbf{u}^{\otimes m}$, and dimension
\begin{align}
\Delta:={n}+{n}^2+\dots {n}^{N}=\frac{{n}^{N+1}-{n}}{{n}-1}.\label{eq:delta_def}
\end{align}
Then the evolution of $\mathbf{z}_i$ is analogous to that of $\mathbf{y}_i$ in \eqref{eq:dxidt_1}, apart from the truncation:
\begin{align}
\frac{{\rm d}\mathbf{z}_i}{{\rm d}t}=
\begin{cases}
    {A}_{0}^1+{A}_{1}^1\mathbf{z}_1+{A}_{2}^1\mathbf{z}_2 &\quad\:\:\:\text{    for } i=1,\\
    {A}_{i-1}^i\mathbf{z}_{i-1}+{A}_{i}^i\mathbf{z}_i+{A}_{i+1}^i\mathbf{z}_{i+1}  &\quad\:\:\:\text{    for } 1<i<N,\\
    {A}_{i-1}^i\mathbf{z}_{i-1}+{A}_{i}^i\mathbf{z}_i &\quad\:\:\:\text{    for } i=N.
\end{cases}
\label{eq:dzidt}
\end{align}  
We now may write \eqref{eq:dzidt} as
\begin{align}
\frac{{\rm d}\mathbf{z}}{{\rm d}t}={\mathcal{A}}_N\mathbf{z}+{b}_N,\quad \quad 
\mathbf{z}(0)=
[\mathbf{u}_0^\mathsf{T},\dots , (\mathbf{u}_0^{\otimes N})^\mathsf{T}]^\mathsf{T},\label{eq:carlemanlinearizedeq}
\end{align}
or explicitly
\begin{align}
\frac{{\rm d}}{{\rm d}t}\left(\begin{array}{c}
    \mathbf{z}_1 \\
    \mathbf{z}_2 \\
    \mathbf{z}_3\\
    \vdots \\
    \mathbf{z}_{N-1}\\[5pt]
    \mathbf{z}_{N}
    \end{array}\right) = 
    \begin{pmatrix}
    {A}_1^1  &{A}_2^1 & & & & \\
    {A}_1^2  &{A}_2^2 &{A}_3^2&&&\\
    &{A}_2^3 &{A}_3^3 &{A}_4^3&&\\
     & &\ddots&\ddots&\ddots&\\
    & & & {A}_{N-2}^{N-1}&{A}_{N-1}^{N-1} &{A}_{N}^{N-1}\\[5pt]
    & & & &{A}_{N}^{N-1} &{A}_{N}^{N}
    \end{pmatrix}
    \left(\begin{array}{c}
    \mathbf{z}_1 \\
    \mathbf{z}_2 \\
    \mathbf{z}_3\\
    \vdots \\
    \mathbf{z}_{N-1}\\[5pt]
    \mathbf{z}_{N}
    \end{array}\right)+
    \left(\begin{array}{c}
    {F}_0 \\
    0 \\
    0\\
    \vdots \\
    0\\[5pt]
    0
    \end{array}\right).
    \label{eq:carlemanmatrix_1}
    \end{align}
Note that ${A}_0^1={F}_0$ in ${b}_N$ above. This is now a finite dimensional linear ODE system, which has an analytical solution given by
\begin{align}
\mathbf{z}(t)=\exp({\mathcal{A}}_N t)\mathbf{z}_0+\left(\int_0^t\exp({\mathcal{A}}_Ns){\rm d}s\right){b}_N.
    \label{eq:analytical_linear_sol}
\end{align}
In practical use cases of Carleman linearization time is discretized, and the above formal solution is approximated by expanding the exponentials up to a certain order, in order to evolve $\mathbf{z}$ between time-steps consecutively.

\section{Main result: Carleman linearization of PDEs}\label{sec:Carleman_for_pdes}
\subsection{Formulation}\label{sec:formulation_pde}
In this section we present our main results. The fundamental idea behind extending Carleman linearization to PDEs is to make use of the continuous analogue of the Kronecker product. For scalar valued functions that is simply the product of the function with itself, with each term having a different copy of the independent (coordinate) variable. Generally, for vector valued functions, i.e. vector fields, it is the Kronecker product of the field with itself, with each term again having a different copy of the independent variable.

The most general quadratic PDE that evolves an $n$ dimensional vector valued function in an $m$ dimensional coordinate space, i.e. $\mathbf{u}(\mathbf{x},t): \mathbb{R}^{m}\times \mathbb{R}_{\geq 0}\to \mathbb{R}^n$, can be written as
\begin{align}
\frac{\partial \mathbf{u}(\mathbf{x})}{\partial t}= F_0(\mathbf{x})+ F_1(\mathbf{x})\mathbf{u}(\mathbf{x}) + F_2(\mathbf{x};\mathbf{w})
\left[\mathbf{u}(\mathbf{x}) \otimes \mathbf{u}(\mathbf{w})\right],\quad
\mathbf{u}(\mathbf{x},t=0)=\mathbf{u}_0(\mathbf{x})
.\label{eq:du(x)_dt_1}
\end{align}
We do not write out the explicit time dependence of the state variables for readability. Again, for simplicity we assume that the matrices, $F_j \in\mathbb{R}^{n\times n^j}$ with $j=0,\:1,\:2$, are time independent. The $\mathbf{x}$ and $\mathbf{w}$ dependencies mean that each of the entries of $F_1(\textbf{x})$ and $F_2(\textbf{x};\textbf{w})$ is a linear partial differential operator with respect to these variables. We use ";" to distinguish dependence on free and eliminated-, i.e. contracted, variables. Above, $\mathbf{w}\in \mathbb{R}^{m}$ is a dummy variable, which can be though of being a copy of $\mathbf{x}$, while $\mathbf{x}$ is a free one. The matrix (vector to be exact) $F_0(\textbf{x})$ is $\textbf{x}$ dependent but contains no differential operators.

Similarly to how it is done in the regular Carleman linearization, we define the auxiliary variables
\begin{align}
\mathbf{y}_i:= \mathbf{u}(\mathbf{x}_1) \otimes \mathbf{u}(\mathbf{x}_2)\otimes\dots \otimes \mathbf{u}(\mathbf{x}_i).\label{eq:y_i_def_pde}
\end{align}
This gives the conceptual basis of the whole generalization procedure. Hence $\mathbf{y}_i: \mathbb{R}^{i\cdot m}\times \mathbb{R}\to \mathbb{R}^{n^{i}}$ is a function of the first $i$ (free) copies of $\mathbf{x}$, which we collectively denote by the vector $\mathbfcal{X}_i:=[\mathbf{x}_1^\mathsf{T},\dots,\mathbf{x}_i^\mathsf{T}]^\mathsf{T}\in \mathbb{R}^{i\cdot m}$ .

We also need to define the variable change operator $\overrightarrow{|}_{(\cdot) \:= \:(\cdot)}$ which may be done through its action on a test function $\mathbf{f}: \mathbb{R}^{m} \to \mathbb{R}^n$ as
\begin{align}
\overrightarrow{\big|}_{\mathbf{a}=\mathbf{b}}\:\mathbf{f}(\mathbf{a}):=\mathbf{f}(\mathbf{b}),
\end{align}
where $\mathbf{a},\mathbf{b}\in\mathbb{R}^{m} $ are the independent variables. This is a linear operator that acts to its right, just as differential operators do. It replaces its first variable within objects to its right with its second one. We can also represent the variable change operator using Dirac deltas as
\begin{align}
\overrightarrow{\big|}_{\mathbf{a}=\mathbf{b}}\:\mathbf{f}(\mathbf{a})= \int_{\mathbb{R}^{m}} \delta^{(m)}(\mathbf{a}-\mathbf{b})\mathbf{f}(\mathbf{a}){\rm d}\mathbf{a}.
\end{align}

\begin{proposition} We claim that the time evolution of $\mathbf{y}_i$, defined in \eqref{eq:y_i_def_pde}, is given by the tridiagonal linear PDE system
\begin{align}
\frac{\partial \mathbf{y}_i}{\partial t}=A_{i-1}^i(\mathbfcal{X}_{i})\mathbf{y}_{i-1}+A_{i}^i(\mathbfcal{X}_i)\mathbf{y}_{i} +A_{i+1}^i(\mathbfcal{X}_{i};\mathbf{x}_{i+1})\mathbf{y}_{i+1},\label{eq:tridiagonal_pde_inf}
\end{align}
with linear operators
\begin{subequations}
\begin{align}
A_{i-1}^i(\mathbfcal{X}_{i}) &:= \sum_{\nu=1}^{i} \Big(\mathbb{I}_{n } \otimes \dots \otimes \underbrace{F_0(\mathbf{x}_\nu)}_{\nu-\text {th position }} \otimes \dots \otimes \mathbb{I}_{n }\Big)
\overrightarrow{\Big|}_{\mathbf{x}_{j}=\mathbf{x}_{j+1}\,\forall\, j\geq \nu}\label{eq:A_i_i_1}\\
A_{i}^i(\mathbfcal{X}_i)&:=\sum_{\nu=1}^{i} \Big(\mathbb{I}_{n } \otimes \dots \otimes \underbrace{F_1(\mathbf{x}_\nu)}_{\nu-\text {th position }} \otimes \dots \otimes \mathbb{I}_{n }\Big)\label{eq:A_i_i}\\
A_{i+1}^i(\mathbfcal{X}_{i};\mathbf{x}_{i+1})&:=\sum_{\nu=1}^{i} \Big(\mathbb{I}_{n } \otimes \dots \otimes \underbrace{F_2(\mathbf{x}_\nu;\mathbf{w}_i)}_{\nu-\text {th position }} \otimes \dots \otimes \mathbb{I}_{n }\Big)
\overrightarrow{\Big|}_{\mathbf{x}_j=\mathbf{x}_{j-1} \,\forall\, j\geq \nu+2}
\overrightarrow{\Big|}_{\mathbf{x}_{\nu+1}=\mathbf{w}_i}
,\label{eq:A_i_i__1}
\end{align}
\end{subequations}
which each have $i$ terms in the Kronecker products.

\begin{proof}
Differentiating \eqref{eq:y_i_def_pde}, applying the Leibniz rule and substituting \eqref{eq:du(x)_dt_1}, we get
\begin{align}
\begin{split}
\frac{\partial \mathbf{y}_i}{\partial t}&= \frac{\partial \mathbf{u}(\mathbf{x}_1)}{\partial t} \otimes  \mathbf{u}(\mathbf{x}_2)\otimes\dots \otimes \mathbf{u}(\mathbf{x}_i) + \dots + \mathbf{u}(\mathbf{x}_1) \otimes \dots \otimes \mathbf{u}(\mathbf{x}_{i-1})\otimes \frac{\partial \mathbf{u}(\mathbf{x}_i)}{\partial t}\\
&=\sum_{\nu=1}^{i}  \mathbf{u}(\mathbf{x}_1) \otimes \dots \otimes 
\Big(F_0(\mathbf{x}_\nu)+ F_1(\mathbf{x}_\nu)\mathbf{u}(\mathbf{x}_\nu) + F_2(\mathbf{x}_\nu;\mathbf{w}_i)
\left[\mathbf{u}(\mathbf{x}_\nu) \otimes \mathbf{u}(\mathbf{w}_i)\right]
\Big)
\otimes \dots \otimes \mathbf{u}(\mathbf{x}_i).
\end{split}
\end{align}
Factoring out the $\mathbf{u}$'s then gives
\begin{subequations}
\begin{align}
&\frac{\partial \mathbf{y}_i}{\partial t}=\\
&\: \phantom{+}\sum_{\nu=1}^{i} \Big(\mathbb{I}_{n } \otimes \dots \otimes F_0(\mathbf{x}_\nu) \otimes \dots \otimes \mathbb{I}_{n }\Big)\Big(\mathbf{u}(\mathbf{x}_1)\otimes \dots \otimes \mathbf{u}(\mathbf{x}_{\nu-1})\otimes 1\otimes \mathbf{u}(\mathbf{x}_{\nu+1})\otimes \dots \otimes\mathbf{u}(\mathbf{x}_i)  \Big)\label{eq:F0_line}\\
&\: +\sum_{\nu=1}^{i} \Big(\mathbb{I}_{n } \otimes \dots \otimes F_1(\mathbf{x}_\nu) \otimes \dots \otimes \mathbb{I}_{n }\Big)\Big(\mathbf{u}(\mathbf{x}_1)\otimes \dots \otimes\mathbf{u}(\mathbf{x}_i)  \Big)\label{eq:F1_line}\\
&\: +\sum_{\nu=1}^{i} \Big(\mathbb{I}_{n } \otimes \dots \otimes F_2(\mathbf{x}_\nu
;\mathbf{w}_i) \otimes \dots \otimes \mathbb{I}_{n }\Big)\Big(\mathbf{u}(\mathbf{x}_1)\otimes \dots \otimes \mathbf{u}(\mathbf{x}_{\nu})\otimes \mathbf{u}(\mathbf{w}_i)\otimes \mathbf{u}(\mathbf{x}_{\nu+1})\otimes \dots \otimes\mathbf{u}(\mathbf{x}_i)  \Big).\label{eq:F2_line}
\end{align}
\end{subequations}
In order to fully factor out the operators, including the summations over $\nu$ (that specify which $\mathbf{u}$ is being evolved in the Kronecker products), one needs to rewrite the right terms in the matrix multiplications in the above equations as $\mathbf{y}_{i-1},\mathbf{y}_{i},\mathbf{y}_{i+1}$, respectively, acted on by variable change operators. In \eqref{eq:F0_line}, the $1$ in the Kronecker product can be ignored, and the resulting sequence is equal to $\mathbf{y}_{i-1}$ after a variable change of $\overrightarrow{|}_{\mathbf{x}_{j}=\mathbf{x}_{j+1}\,\forall\, j\geq \nu}$. In \eqref{eq:F1_line} the right term in the product is already $\mathbf{y}_i$ so no change has to be made. In \eqref{eq:F2_line} the right term is equal to $\mathbf{y}_{i+1}$ after the variable changes $\overrightarrow{|}_{\mathbf{x}_j=\mathbf{x}_{j-1} \,\forall\, j\geq \nu+2}
\overrightarrow{|}_{\mathbf{x}_{\nu+1}=\mathbf{w}_i}$ (with the second operator acting first). Then we may write
\begin{align}
\begin{split}
\frac{\partial \mathbf{y}_i}{\partial t}
&=\sum_{\nu=1}^{i} \Big(\mathbb{I}_{n } \otimes \dots \otimes F_0(\mathbf{x}_\nu) \otimes \dots \otimes \mathbb{I}_{n }\Big)
\overrightarrow{\Big|}_{\mathbf{x}_{j}=\mathbf{x}_{j+1}\,\forall\, j\geq \nu}\:
\Big(\mathbf{u}(\mathbf{x}_1)\otimes \dots \otimes\mathbf{u}(\mathbf{x}_{i-1})  \Big)\\
&+\sum_{\nu=1}^{i} \Big(\mathbb{I}_{n } \otimes \dots \otimes F_1(\mathbf{x}_\nu) \otimes \dots \otimes \mathbb{I}_{n }\Big)
\Big(\mathbf{u}(\mathbf{x}_1)\otimes \dots \otimes\mathbf{u}(\mathbf{x}_i)  \Big)\\
&+\sum_{\nu=1}^{i} \Big(\mathbb{I}_{n } \otimes \dots \otimes F_2(\mathbf{x}_\nu
;\mathbf{w}_i) \otimes \dots \otimes \mathbb{I}_{n }\Big)
\overrightarrow{\Big|}_{\mathbf{x}_j=\mathbf{x}_{j-1} \,\forall\, j\geq \nu+2}
\overrightarrow{\Big|}_{\mathbf{x}_{\nu+1}=\mathbf{w}_i}\:
\Big(\mathbf{u}(\mathbf{x}_1)\otimes \dots \otimes\mathbf{u}(\mathbf{x}_{i+1})  \Big)\\[5pt]
&=A_{i-1}^i(\mathbfcal{X}_{i})\mathbf{y}_{i-1}+A_{i}^i(\mathbfcal{X}_i)\mathbf{y}_{i} +A_{i+1}^i(\mathbfcal{X}_{i};\mathbf{x}_{i+1})\mathbf{y}_{i+1}
,\label{eq:proof_line_2}
\end{split}
\end{align}
making use of the matrix definitions from (\ref{eq:A_i_i_1}-\ref{eq:A_i_i__1}) to get the last line. To be explicit, we have $\mathbf{x}_1,\dots,\mathbf{x}_i$ as free variables in the above, while $\mathbf{x}_{i+1}$ is eliminated by the variable change operator in $A_{i+1}^i(\mathbfcal{X}_{i};\mathbf{x}_{i+1})$. In turn, $A_{i+1}^i(\mathbfcal{X}_{i};\mathbf{x}_{i+1})$ introduces the dummy variable $\mathbf{w}_i$ which is eliminated through contraction by $F_2(\cdot;\mathbf{w}_i)$. Note that the index $i$ on $\mathbf{w}_i$ becomes essential when exponentiating the matrix arising from this linear system, because that exponential generally contains products of $A_{i+1}^i(\mathbfcal{X}_{i};\mathbf{x}_{i+1})$ matrices with different $i$'s. This way there is no ambiguity among the different $\mathbf{w}_i$'s.
\end{proof}
\end{proposition}

\subsection{Truncation}\label{sec:truncation_pde}
We truncate \eqref{eq:tridiagonal_pde_inf} at level $N$ to make the resulting system solvable in practise. We define the finite state vector $\mathbf{z}(\mathbfcal{X}_N):=[\mathbf{z}_1^\mathsf{T},\dots , \mathbf{z}_N^\mathsf{T}]^\mathsf{T}:  \mathbb{R}^{m\cdot n}\times \mathbb{R}_{\geq 0} \to \mathbb{R}^{\Delta}$ that has block elements $\mathbf{z}_i:=\mathbf{u}(\mathbf{x}_1)\otimes \dots \otimes \mathbf{u}(\mathbf{x}_i)$, and dimension $\Delta$ from \eqref{eq:delta_def}. It is evolved by the truncated system
\begin{align}
\frac{\partial \mathbf{z}_i}{\partial t}
=
\begin{cases}
    {A}_{0}^1(\mathbfcal{X}_{1})+{A}_{1}^1(\mathbfcal{X}_{1})\mathbf{z}_1+{A}_{2}^1(\mathbfcal{X}_{1};\mathbf{x}_2)\mathbf{z}_2 &\quad\:\:\:\text{    for } i=1,\\
    {A}_{i-1}^i(\mathbfcal{X}_{i})\mathbf{z}_{i-1}+{A}_{i}^i(\mathbfcal{X}_{i})\mathbf{z}_i+{A}_{i+1}^i(\mathbfcal{X}_{i};\mathbf{x}_{i+1})\mathbf{z}_{i+1}  &\quad\:\:\:\text{    for } 1<i<N,\\
    {A}_{N-1}^N(\mathbfcal{X}_{N})\mathbf{z}_{N-1}+{A}_{N}^N(\mathbfcal{X}_{N})\mathbf{z}_N &\quad\:\:\:\text{    for } i=N.
\end{cases}
\label{eq:dzidt_3}
\end{align}
This can be represented as a linear PDE system of the form
\begin{align}
\begin{split}
\frac{\partial \mathbf{z}(\mathbfcal{X}_N)}{\partial t}&={\mathcal{A}}_N(\mathbfcal{X}_N)\mathbf{z}(\mathbfcal{X}_N)+{b}_N(\mathbfcal{X}_1),\\[5pt]
\mathbf{z}(\mathbfcal{X}_N,0)&=\mathbf{z}_0(\mathcal{X}_N):=
[\mathbf{u}_0(\mathbf{x}_1)^\mathsf{T},\dots, (\mathbf{u}_0(\mathbf{x}_1)\otimes \dots \otimes \mathbf{u}_0(\mathbf{x}_N))^\mathsf{T} ]^\mathsf{T}
,\label{eq:carlemanlinearizedeq_pde}
\end{split}
\end{align}
with
\begin{align}
\begin{split}
{\mathcal{A}}_N(\mathbfcal{X}_N)&:=
\begin{pmatrix}
    {A}_1^1(\mathbfcal{X}_1)  &{A}_2^1(\mathbfcal{X}_1;\mathbf{x}_2) & & & & \\[5pt]
    {A}_1^2(\mathbfcal{X}_2)  &{A}_2^2(\mathbfcal{X}_2) &{A}_3^2(\mathbfcal{X}_2;\mathbf{x}_3)&&&\\[5pt]
     & \ddots&\ddots&\ddots&&\\[5pt]
    & & & {A}_{N-2}^{N-1}(\mathbfcal{X}_{N-1})&{A}_{N-1}^{N-1}(\mathbfcal{X}_{N-1}) &{A}_{N}^{N-1}(\mathbfcal{X}_{N-1};\mathbf{x}_N)\\[5pt]
    & & & &{A}_{N-1}^{N}(\mathbfcal{X}_N) &{A}_{N}^{N}(\mathbfcal{X}_N)
    \end{pmatrix}\\
{b}_N(\mathbfcal{X}_1)&:=
 \left(\begin{array}{c}
    {F}_0(\mathbfcal{X}_1) \\
    0 \\
    \vdots \\
    0\\[5pt]
    0
    \end{array}\right).
    \label{eq:carlemanmatrix_2}
\end{split}
\end{align}
This is now a finite dimensional, and has an analytical solution given by
\begin{align}
\begin{split}
\mathbf{z}(\mathbfcal{X}_N,t)=\exp\Big({\mathcal{A}}_N(\mathbfcal{X}_N) t\Big)\mathbf{z}_0(\mathbfcal{X}_N)+\left(\int_0^t\exp\Big({\mathcal{A}}_N(\mathbfcal{X}_N)s\Big){\rm d}s\right){b}_N(\mathbfcal{X}_1).
\label{eq:analytical_linear_sol_pde}
\end{split}
\end{align}

\section{Application to toy problems}\label{sec:toy_problems}
\subsection{Example: Burger's equation}
Burger's equation is a model of fluids. In one dimension it evolves the scalar field $u(x,t)$ according to
\begin{align}
\frac{\partial u}{\partial t}=\mu \frac{\partial^2 u}{\partial x^2}-u\frac{\partial u}{\partial x},\quad\quad u(x,t=0)=u_0(x)\label{eq:burger_def}
\end{align}
where $x\in \mathbb{R}$, $t\geq 0$ and $\mu\geq 0$ is a constant. Since $u$ is a scalar, the $F_j$ and corresponding the $A_k^i$ matrices are also one dimensional. From the above equation, we identify
\begin{align}
F_2(x;w)=- \overrightarrow{\Big|}_{w=x}\frac{\partial}{\partial x},\quad F_1(x)=\mu \frac{\partial^2}{\partial x^2},\quad F_0(x)=0.\label{eq:burger_F_j}
\end{align}
This way we have
\begin{align}
\begin{split}
F_2(x;w) u(x)u(w)&= -\overrightarrow{\Big|}_{w=x}\frac{\partial}{\partial x}u(x)u(w)\\
&=-\overrightarrow{\Big|}_{w=x} u(w)\frac{\partial u(x)}{\partial x}\\
&= -u(x)\frac{\partial u(x)}{\partial x}
\end{split}
\end{align}
Note that
\begin{align}
F_2(x;w) = -\overrightarrow{\Big|}_{w=x}\frac{\partial }{\partial w}
\end{align}
gives the same result. The $F_2$ operator is not unique generally, due to the structure of the Kronecker power. Then, substituting \eqref{eq:burger_F_j} into (\ref{eq:A_i_i_1}-\ref{eq:A_i_i__1}), we obtain
\begin{align}
\begin{split}
A_{i-1}^i(\mathbfcal{X}_{i}) &=0\\
A_{i}^i(\mathbfcal{X}_i)&=\sum_{\nu=1}^{i} F_1(x_\nu) = \mu\sum_{\nu=1}^i  \frac{\partial^2}{\partial x_\nu ^2}\\
A_{i+1}^i(\mathbfcal{X}_{i};x_{i+1})&=\sum_{\nu=1}^{i} F_2(x_\nu;w_i)
\overrightarrow{\Big|}_{x_j=x_{j-1} \,\forall\, j\geq \nu+2}
\overrightarrow{\Big|}_{x_{\nu+1}=w_i}
=  - \sum_{\nu=1}^{i} \overrightarrow{\Big|}_{w_i=x_\nu}\frac{\partial}{\partial x_\nu}
\overrightarrow{\Big|}_{x_j=x_{j-1} \,\forall\, j\geq \nu+2}
\overrightarrow{\Big|}_{x_{\nu+1}=w_i}
,\label{eq:A_i_i_burger}
\end{split}
\end{align}
from which we can construct the linear system given by \eqref{eq:carlemanmatrix_2}. The exponential of $\mathcal{A}_N(\mathbfcal{X}_N)t$ cannot be evaluated analytically, but the system can be numerically evolved with various time integrators.

For the special case of $\mu=0$ (inviscid Burger's equation), we have $F_1(x)=0$ and it is possible to express the formal solution \eqref{eq:analytical_linear_sol_pde} directly. We also have $b_N=0$. The matrix ${\mathcal{A}}_N(\mathbfcal{X}_N)$ takes the form
\begin{align}
{\mathcal{A}}_N(\mathbfcal{X}_N)=
\begin{pmatrix}
0 & A_{2}^1(\mathbfcal{X}_{1};x_{2}) & & &\\
 & 0&A_{3}^2(\mathbfcal{X}_{2};x_{3}) & &\\
&&\ddots &\ddots &\\
&& &0&A_{N}^{N-1}(\mathbfcal{X}_{N-1};x_{N})\\
&&&&0
\end{pmatrix}.
\end{align}
Exponentiating, we get
\begin{align}
\exp\left({\mathcal{A}}_N(\mathbfcal{X}_N)t\right) = 
\begin{pmatrix}
1\:\: & \dfrac{t}{1!}A_{2}^1 \:\:&\dfrac{t^2}{2!}A_{2}^1 A_{3}^2 \:\:&\dots  \:\:&\dfrac{t^{N-1}}{(N-1)!}A_{2}^1 \dots A_{N}^{N-1} \\[8pt]
 & 1&\dfrac{t}{1!} A_{3}^2 &\dots &\dfrac{t^{N-2}}{(N-2)!}A_{3}^2\dots A_{N}^{N-1}\\[8pt]
&&\ddots &\ddots &\vdots\\[8pt]
 &  &&1 &\dfrac{t}{1!} A_{N}^{N-1} \\[8pt]
& & &  & 1 \\[8pt]
\end{pmatrix}
\end{align}
where the coordinate dependence of the operators was left out for readability. The elements read 
\begin{align}
\Big (\exp\left({\mathcal{A}}_N(\mathbfcal{X})t\right) \Big)_{ij} = 
\begin{cases}
    0&\quad\:\:\:\text{for  } i>j,\\
    1  &\quad\:\:\:\text{for  } i=j,\\
    \dfrac{t^{j-i}}{(j-i)!}A_{i+1}^i(\mathbfcal{X}_{i};x_{i+1})
    \dots A_{j}^{j-1}(\mathbfcal{X}_{j-1};x_{j})&\quad\:\:\:\text{for  } i<j.\label{eq:ddv}
\end{cases}
\end{align}
Therefore we can formally write the solution of \eqref{eq:burger_def} as
\begin{align}
u(x_1,t)&=\sum_{j=1}^N  \Big (\exp\left({\mathcal{A}}_N(\mathbfcal{X})t\right) \Big)_{1,j}  \Big(u_0(x_1)\dots u_0(x_j)\Big)\\
&= u_0(x_1) + \sum_{j=2}^N \dfrac{t^{j-1}}{(j-1)!}
A_{2}^1(\mathbfcal{X}_{1};x_{2}) \dots A_{j}^{j-1}(\mathbfcal{X}_{j-1};x_{j})\Big(u_0(x_1)\dots u_0(x_j)\Big).
\end{align}
This takes the form of Taylor expansion in $t$ and may form the basis of a time-integrator.

\subsection{Example: Vlasov equation}
The Vlasov equation is a widely used model of plasmas that evolves the phase-space distribution function of each particle species of the physical system. We focus on the one dimensional model, i.e. we have one dimension of position $x^{(1)}\in \mathbb{R}$ and one of velocity $x^{(2)}\in \mathbb{R}$ making up the independent variable $\mathbf{x}=[x^{(1)},x^{(2)}]^\mathsf{T}\in \mathbb{R}^2$. We furthermore choose to work with 2 species, i.e. the dependent variable is $\mathbf{u}(\mathbf{x})=[u_1(\mathbf{x}),u_2(\mathbf{x})]^\mathsf{T}\in \mathbb{R}^2_{\geq 0}$, with $u_1(\textbf{x})\in \mathbb{R}_{\geq 0}$ and $u_2(\textbf{x})\in \mathbb{R}_{\geq 0}$ being the individual distribution functions. Using Gauss's law to express the electric field that accelerates the plasma, we can formulate the system as
\begin{align}
\begin{split}
\frac{\partial u_1(\mathbf{x})}{\partial t}&=-x^{(2)}\frac{\partial u_1(\mathbf{x})}{\partial x^{(1)}} +c_1 \frac{\partial u_1(\mathbf{x})}{\partial x^{(2)}}\int^{x^{(1)}}_{-\infty }{\rm d}w^{(1)} \int^{\infty}_{-\infty }{\rm d}w^{(2)}
\big[u_2(\mathbf{w})-u_1(\mathbf{w})\big]\\
\frac{\partial u_2(\mathbf{x})}{\partial t}&=-x^{(2)}\frac{\partial u_2(\mathbf{x})}{\partial x^{(1)}} +c_2 \frac{\partial u_2(\mathbf{x})}{\partial x^{(2)}}\int^{x^{(1)}}_{-\infty }{\rm d}w^{(1)} \int^{\infty}_{-\infty }{\rm d}w^{(2)}
\big[u_2(\mathbf{w})-u_1(\mathbf{w})\big],\label{eq:vlasov}
\end{split}
\end{align}
where $c_1,c_2\in \mathbb{R}$ are species dependent constants. Note that the integrals over $w^{(1)}$ are cumulative. Then, we read off that
\begin{align}
F_0(\mathbf{x})=
\begin{pmatrix}
    0\\
    0
\end{pmatrix}
, \quad\quad F_1(\mathbf{x})=
\begin{pmatrix}
   -x^{(2)}\dfrac{\partial }{\partial x^{(1)}} &0  \\
    0& -x^{(2)}\dfrac{\partial }{\partial x^{(1)}}
\end{pmatrix}= - \mathbb{I}_2\: x^{(2)}\frac{\partial }{\partial x^{(1)}} .
\end{align}
$F_2(\mathbf{x};\mathbf{w})$ can be found by manipulating the nonlinear part of \eqref{eq:vlasov} into the form given in \eqref{eq:du(x)_dt_1} as
\begin{align}
\begin{split}
F_2(\mathbf{x};\mathbf{w})
\left[\mathbf{u}(\mathbf{x}) \otimes \mathbf{u}(\mathbf{w})\right] &=
\begin{pmatrix}
c_1 \dfrac{\partial u_1(\mathbf{x})}{\partial x^{(2)}}\int^{x^{(1)}}_{-\infty }{\rm d}w^{(1)} \int^{\infty}_{-\infty }{\rm d}w^{(2)}
\big[u_2(\mathbf{w})-u_1(\mathbf{w})\big] \\[8pt]
c_2 \dfrac{\partial u_2(\mathbf{x})}{\partial x^{(2)}}\int^{x^{(1)}}_{-\infty }{\rm d}w^{(1)} \int^{\infty}_{-\infty }{\rm d}w^{(2)}
\big[u_2(\mathbf{w})-u_1(\mathbf{w})\big]
\end{pmatrix}\\
&=
\begin{pmatrix}
\left[F_2(\mathbf{x};\mathbf{w})\right]_{1,1}&\left[F_2(\mathbf{x};\mathbf{w})\right]_{1,2}&0&0\\[5pt]
0&0&\left[F_2(\mathbf{x};\mathbf{w})\right]_{2,3}&\left[F_2(\mathbf{x};\mathbf{w})\right]_{2,4}
\end{pmatrix}
\begin{pmatrix}
u_1(\mathbf{x}) u_1(\mathbf{w})\\
u_1(\mathbf{x}) u_2(\mathbf{w})\\
u_2(\mathbf{x}) u_1(\mathbf{w})\\
u_2(\mathbf{x}) u_2(\mathbf{w})\\
\end{pmatrix},
\end{split}
\end{align}
with matrix entries
\begin{align}
\begin{split}
\left[F_2(\mathbf{x};\mathbf{w})\right]_{1,1} &=- c_1 \dfrac{\partial }{\partial x^{(2)}}\int^{x^{(1)}}_{-\infty }{\rm d}w^{(1)} \int^{\infty}_{-\infty }{\rm d}w^{(2)}\\
\left[F_2(\mathbf{x};\mathbf{w})\right]_{1,2} &= - \left[F_2(\mathbf{x};\mathbf{w})\right]_{1,1}\\
\left[F_2(\mathbf{x};\mathbf{w})\right]_{2,3} &=-c_2 \dfrac{\partial} {\partial x^{(2)}}\int^{x^{(1)}}_{-\infty }{\rm d}w^{(1)} \int^{\infty}_{-\infty }{\rm d}w^{(2)}\\
\left[F_2(\mathbf{x};\mathbf{w})\right]_{2,4} &=-\left[F_2(\mathbf{x};\mathbf{w})\right]_{2,3}.
\end{split}
\end{align}
These entries all eliminate $\mathbf{w}$ but keep $\mathbf{x}$ free as required.

\section{Discussion and conclusion}\label{sec:disc_conlc}
In Section \ref{sec:formulation_pde} we presented our generalization of Carleman linearization, that embeds systems of quadratically nonlinear PDEs of the form of \eqref{eq:du(x)_dt_1}, into infinite dimensional linear PDEs of the form of \eqref{eq:tridiagonal_pde_inf}. In Section \ref{sec:truncation_pde}, we showed how the linear system can be truncated, making the resulting equation solvable on a computer.

When solving nonlinear PDE systems algorithmically, our method has the advantage over finite-difference-based ones that it does not require the discretization of coordinate space, hence eliminating an error that is typically the largest in conventional methods for these problems in numerical analysis. Although an error analysis has not yet been performed, the drawback of our method is expected to be a strong requirement on the smoothness of the solution and the initial condition, as well as on the form of the partial differential operators involved. We anticipate that our method might be useful in the future for simulating certain nonlinear PDE systems on analog quantum devices with a potential computational advantage over classical and digital-quantum approaches.

The next step to be done is performing a proper analysis of the error arising from truncating the Carleman linearized system, similar to what has been done for the original Carleman linearization in literature. The error bound on the numerical solution at time $T$ presented in Ref. \cite{forets2017explicit}, depends on the quantities $T,\| F_2\|,\mu(F_1),\|\mathbf{u}_t\|$ and $N$, where $\mu(\cdot)$ is the log-norm and $\|\mathbf{u}_t\|$ is the maximum of the norm of the state vector over $t\in[0,T]$. The conclusion is that if (a): the system is overall dissipative (i.e. $\mu(F_1)<0$); and (b) $\| F_2\|+\| F_0\|$ is small compared to $|\mu(F_1)|$; then the error goes to zero as $N\to \infty$. We expect a qualitatively similar error analysis to be possible in the PDE case, i.e. formulating an error bound using the continuous analogues of the mentioned quantities.

Finally, we note that the framework of Carleman linearization also provides a procedural approach to performing quadratization on higher-degree polynomial nonlinearities, as shown in \cite[Section 3.3]{forets2017explicit}. This method forms the basis of the quantum algorithm proposed in \cite{surana2023efficient}. Quadratization is a well-known technique for ODEs (primarily without Carleman linearization) \cite{quadratization1,quadratization2,quadratization3,quadratization4}, however it has not yet been extended to PDEs to the best of the authors' knowledge. We anticipate that our framework may also be effective for this purpose, though this remains for future investigation.

\bibliographystyle{IEEEtran}
\bibliography{ref}

\end{document}